\documentclass[11pt,reqno]{amsart}
\usepackage{amssymb}
\usepackage[all]{xy}
\newtheorem{thm}{Theorem}[section]
\newtheorem{lem}[thm]{Lemma}
\newtheorem{cor}[thm]{Corollary}
\newtheorem{con}[thm]{Conjecture}
\begin{document}
\title[Bordism classes of the multiple points manifolds]
{Bordism classes of the multiple points manifolds\\ of smooth immersion}
\author{Konstantin Salikhov}
\thanks{Partially supported by the Russian Fond for Basic Research
Grant No. 99-01-00009}
\address{Dept. of Differential Geometry, Faculty of Mechanics and
Mathematics, Moscow State University, Moscow, 119899, Russia}
\email{salikhov@mccme.ru}
\keywords{Immersion, multiple points of immersion, bordism group of immersion}
\subjclass{Primary: 57R42}
\begin{abstract}
Let $f:V^n\looparrowright M^m$ be a smooth generic immersion. Then the set
of points, that have at least $k$ preimages is an image of a (non-generic)
immersion. If the manifolds $V^n$ and $M^m$ are oriented and $m-n$ is even,
then the manifold of $k$-fold points is also oriented. In this paper we
compute the oriented bordism class of the manifold of $k$-fold points
in terms of the differential $df$, provided the tangent
bundle of the manifold $M^m$ has a nowhere zero cross-section.
\end{abstract}
\maketitle
\section{Introduction}
In this paper we will consider smooth orientable $C^\infty$-manifolds and
$C^\infty$-mappings between them. Let $V^n$ and $M^m$ be manifolds without
boundary, $V^n$ be compact, $f:V^n\looparrowright M^m$ be a smooth generic
immersion and $m-n$ be {\it even}. From Thom multijet transversality
theorem \cite{GG} it follows that the set of immersions $f$ such that the
map $f^{(k)}:V^{(k)}\to M^{(k)}$ is transverse to the "thin" diagonal
$\Delta_k(M)=\{(x,\dots,x)\mid x\in M\}\subset M^{(k)}$
outside the "thick" diagonal
$\Delta_2(V)=\{(x_1,\dots,x_k)\mid\exists\ i\neq j: x_i=x_j\}\subset V^{(k)}$
is open and everywhere dense in the set of all immersions
$V^n\looparrowright M^m$ (in $C^\infty$ Whitney topology). Therefore,
for a generic $f$ the set
$V_k=(f^{(k)})^{-1}(\Delta_k(M))\backslash\Delta_2(V)\subset V^{(k)}$
is an oriented submanifold. Denote by $\Sigma_k$ the group of permutations
on $k$ elements. Permutation of factors in the product $V^{(k)}$ induce
the free action of $\Sigma_k$ on the set $V^{(k)}$. Denote by
$\widetilde{V}_k$ and $\widetilde{M}_k$ the quotient manifolds
$V_k/\Sigma_{k-1}$ and $V_k/\Sigma_k$, where the subgroup
$\Sigma_{k-1}\subset\Sigma_k$ is the stabilizer of the first element.
Since $m-n$ is even, the action of $\Sigma_k$ preserves the orientation
of the manifold $V_k$. Therefore $\widetilde{V}_k$ and
$\widetilde{M}_k$ are canonically oriented.
Let us define immersions $f_k:\widetilde{V}_k\looparrowright V^n$ and
$g_k:\widetilde{M}_k\looparrowright M^m$ by formulae
$f_k(x_1,[x_2,\dots,x_k])=x_1$ and $g_k[x_1,\dots,x_k]=f(x_1)$
(see \cite{Ronga} for details).

For an integer $k>0$ and an immersion $f:V^n\looparrowright M^m$
we assign the oriented bordism classes
$(\widetilde{V}_k,f_k)\in\Omega_{m-k(m-n)}(V^n)$ and
$(\widetilde{M}_k,g_k)\in\Omega_{m-k(m-n)}(M^m)$. From \cite{KS} it
follows that these classes do not change under a regular homotopy of the
immersion $f$. By the fundamental Smale-Hirsch theorem \cite{Hirsch},
the set of regular homotopy classes of immersions $V^n\looparrowright M^m$
is in 1-1 correspondence with the set of linear monomorphism classes of 
the tangent bundles $\tau V\to\tau M$. The class of the immersion
$f:V^n\looparrowright M^m$ corresponds to the class of the differential
$df:\tau V\to\tau M$. Hence, the classes
$(\widetilde{V}_k,f_k)\in\Omega_{m-k(m-n)}(V^n)$ and
$(\widetilde{M}_k,g_k)\in\Omega_{m-k(m-n)}(M^m)$ have to be computable
in terms of the differential $df:\tau V\to\tau M$. In this paper
(see corollary~\ref{maincor}) we compute $(\widetilde{V}_k,f_k)$ and
$(\widetilde{M}_k,g_k)$ up to elements of order $(k-1)!$ and $k!$,
respectively. This weakening is in some sense natural, for the manifolds
$\widetilde{V}_k$ É $\widetilde{M}_k$ were constructed as the images of
$(k-1)!$- and $k!$-fold coverings.

Recall that all the finite order elements of $\Omega_*(pt)$ have the order~2.
Therefore, for $M^m={\mathbb R}^m$ we compute the classes
$(\widetilde{M}_k,g_k)$ up to elements of order~2. It is another approach
to \cite[Theorem~5]{Sz 98}, where all the Pontrjagin numbers of the manifolds
$\widetilde{M}_k$ were computed in terms of the Pontrjagin classes
of the manifold $V^n$ and the integral Euler class of the normal bundle
of the immersion $f$. Lemma~\ref{double} gives a formula to compute
the unoriented bordism class $(\widetilde{V}_2,f_2)_2\in\Re_*(V^n)$,
if we do not require $m-n$ to be even, and $V^n$ and $M^m$ to be oriented.
For up-to-date reviews of results on the bordism classes of self-intersection
manifolds see \cite{Sz 98} for oriented case, and \cite{AE} for unoriented
case.

\section{Formulation of results}
Denote by $S\tau M$ the spherical fibration, associated to the tangent
bundle $\tau M$ and by $Sdf:S\tau V\to S\tau M$ the fiberwise monomorphism
of spherical fibrations, induced by the differential $df:\tau V\to\tau M$.
Since the manifolds $V^n$, $M^m$, $S\tau V$ and $S\tau M$ are oriented
(in usual sense), they are oriented in {\it oriented bordism theory}
\cite{CF,Sw}. Thus, there is the {\it Poincare duality} on these manifolds.
Denote by $Sdf^!:\Omega_*(S\tau M)\to\Omega_*(S\tau V)$ the Gysin
homomorphism, induced by the mapping $Sdf$.

Let us formulate the key~lemma of this paper.
\begin{lem}\label{double}
Let  $M^m$ be an oriented manifold without boundary such that there is
a nowhere zero cross-section of the tangent bundle $\tau M$, or, in other
words, a section $s_M:M^m\to S\tau M$. Then for any generic immersion
$f:V^n\looparrowright M^m$ of compact oriented manifold without boundary
$V^n$ the bordism class $(\widetilde{V}_2,f_2)\in\Omega_{2n-m}(V^n)$ is
\begin{equation}\label{center}
(\widetilde{V}_2,f_2)=(-1)^{m-1}i_*Sdf^!(M^m,s_M),
\end{equation}
where $i$ is the natural projection $S\tau V\to V^n$.
\end{lem}

To formulate our results, it will be convenient to use the
oriented cobordism classes $v_k\in\Omega^{(k-1)(m-n)}(V^n)$ and
$m_k\in\Omega^{k(m-n)}_{comp.}(M^m)$, the Poincare duals to
$(\widetilde{V}_k,f_k)$ and $(\widetilde{M}_k,g_k)$, respectively.
Denote by $1_V$ the identity element of the ring $\Omega^*(V^n)$,
and by $f_!$ the Gysin homomorphism, induced by the map $f$.
\begin{cor}\label{ecor} Under the conditions of lemma~\ref{double} and
if $m-n$ is even, the Euler class $e$ of the normal bundle of the
immersion $f$ is
\begin{equation*}
e=f^*f_!(1_V)+(-1)^m\gamma i_*Sdf^!(M^m,s_M),
\end{equation*}
where $\gamma:\Omega_*(V^n)\to\Omega^{n-*}(V^n)$ is the Poincare duality.
\end{cor}

\begin{cor}\label{maincor} Under the conditions of lemma~\ref{double} and
if $m-n$ is even,
\begin{equation*}
\begin{gathered}
(k-1)!\cdot v_k=
\varphi_{k-1}\circ\varphi_{k-2}\circ\cdots\circ\varphi_1(1_V)\\
k!\cdot m_k=
f_!\circ\varphi_{k-1}\circ\varphi_{k-2}\circ\cdots\circ\varphi_1(1_V),
\end{gathered}
\end{equation*}
where $\varphi_k(a)=f^*f_!(a)-k\cdot e\cup a$, and $e$ is the Euler class
of the normal bundle of the immersion $f$ (which was computed in
corollary~\ref{ecor}).
\end{cor}

\section{The bordism group of immersions}
Let us call two oriented immersions $f_0:V_0^n\looparrowright M^m$ and
$f_1:V_1^n\looparrowright M^m$ {\it bordant}, if there exists a
compact oriented manifold with boundary $W^{n+1}$ such that
$\partial W^{n+1}=V_0^n\sqcup(-V_1^n)$, and an immersion
$W^{n+1}\looparrowright M^m\times[0,1]$ such that for a collar
$V_0^n\times[0,\varepsilon)\sqcup(-V_1^n)\times(1-\varepsilon,1]$ of the
boundary $\partial W^{n+1}$ the restrictions
$F|_{V_0^n\times[0,\varepsilon)}=f_0\times{\qopname\relax o{id}}$ and
$F|_{(-V_1^n)\times(1-\varepsilon,1]}=f_1\times{\qopname\relax o{id}}$.
Then the set of equivalence classes of bordant oriented
immersions with disjoint union operation is a group
$Imm_n^{SO}(M^m)$. The group $Imm_n^{SO}(M^m)=\left[M^m,Q MSO(m-n)\right]$,
where $Q X=\varinjlim\Omega^q S^q X$ is the infinite loop space of infinite
suspension and $MSO$ is the Thom spectrum \cite{Wells}.

From the results of paper \cite{KS} it follows that the map, assigning
for any immersion $f:V^n\looparrowright M^m$ the bordism class
$(\widetilde{M}_k,g_k)\in\Omega_{m-k(m-n)}(M^m)$, is a well-defined
homomorphism $\varepsilon_k:Imm_n^{SO}(M^m)\to\Omega_{m-k(m-n)}(M^m)$.
The classes $(\widetilde{M}_k,g_k)$ involve much information about the
class of immersion $[f]\in Imm_n^{SO}(M^m)$. The following theorem was
proved in \cite{Sz 80}, using algebraic technics.
\begin{thm}[{\cite[Corollary 1]{Sz 80}}]\label{szthm}
If $3n+1<2m$ and $m-n$ is even, then the homomorphism
\begin{equation*}
\varepsilon_1\oplus\varepsilon_2:Imm_n^{SO}({\mathbb R}^m)\to
\Omega_n\oplus\Omega_{2n-m}
\end{equation*}
is an isomorphism modulo the class $C_2$ of finite 2-primary groups.
\end{thm}

Our calculations probably clarify the geometric core of theorem~\ref{szthm}.
The matter is that in these very dimensional restrictions ($3n+1<2m$)
any {\it skew} map
\footnote{a skew map $\tau V\to\tau M$ can be understood as the
"fiberwise cone" over a fiber map $h:S\tau V\to S\tau M$ such that
$h(-x)=-h(x)$ in each fiber}
$\tau V\to\tau M$ can be homotoped to a monomorphism
of tangent bundles (see details in \cite{HH}). Our formula
(lemma~\ref{double}) connects the map of spherical fibrations, induced by
the differential $df:\tau V\to \tau M$, and the oriented bordism
class $2(\widetilde{M}_2,g_2)\in\Omega_{2n-m}(M^m)$. These reasoning
was the initial motivation of this paper.

\begin{con}\label{conj}
Let  $M^m$ be an oriented manifold without boundary such that there is
a nowhere zero cross-section of the tangent bundle $\tau M$,
$3n+1<2m$, and $m-n$ be even. Then the homomorphism
\begin{equation*}
\varepsilon_1\oplus\varepsilon_2:Imm_n^{SO}(M^m)\to
\Omega_n(M^m)\oplus\Omega_{2n-m}(M^m)
\end{equation*}
is an isomorphism modulo the class $C_2$ of finite 2-primary groups.
\end{con}

\section{Proofs}
Denote the diagonal $\Delta(V)=\{(x,x)\in V\times V\mid x\in V\}$.
Since $f:V^n\looparrowright M^m$ is an immersion, there exists a small
enough tubular neighborhood $U_V$ of the diagonal $\Delta(V)$ in $V\times V$
such that $f^{(2)}(U_V\backslash\Delta(V))\cap\Delta(M)=\emptyset$.
Note that $\partial(V^{(2)}\backslash U_V)=\partial(U_V)$. Since
$f^{(2)}(U_V\backslash\Delta(V))\cap\Delta(M)=\emptyset$, we get the map
$f^{(2)}:(V^{(2)}\backslash U_V,\partial(V^{(2)}\backslash U_V))\to
(M^{(2)},M^{(2)}\backslash\Delta(M))$. Denote by $U_M$ a tubular
neighborhood of the diagonal $\Delta(M)$ in $M^{(2)}$. Without loss
of generality we may assume that $f^{(2)}(U_V)\subset U_M$.
Denote the inclusion
\begin{equation*}
j:(U_M,U_M\backslash\Delta(M))\hookrightarrow
(M^{(2)},M^{(2)}\backslash\Delta(M))
\end{equation*}
By excision axiom \cite{Sw}, the homomorphisms
$j_*:\Omega_*(U_M,U_M\backslash\Delta(M))\to
\Omega_*(M^{(2)},M^{(2)}\backslash\Delta(M))$ and
$j^*:\Omega^*(M^{(2)},M^{(2)}\backslash\Delta(M))\to
\Omega^*(U_M,U_M\backslash\Delta(M))$ are isomorphisms.
Note that the pair $(U_M,\Delta(M))$ is canonically isomorphic to the
pair $(\tau M,\tau_0M)$ \cite{MS}. Since $M^m$ is oriented, there exists
the {\it Thom class} $t\in\Omega^m(\tau M,\tau M\backslash\tau_0M)$
of the tangent bundle $\tau M$.
\begin{lem}\label{lemma}
The class $(\widetilde{V}_2,f_2)\in\Omega_{2n-m}(V^n)$ for an immersion
$f:V^n\looparrowright M^m$ can be calculated in the following way
\begin{equation}\label{formula}
(\widetilde{V}_2,f_2)=(\pi_1)_*\left((f^{(2)})^*((j^*)^{-1}t)\cap
\left[V^{(2)}\backslash U_V,\partial(V^{(2)}\backslash U_V)\right]\right),
\end{equation}
where $\pi_1:V^{(2)}\backslash U_V\to V$ is the projection on the first
factor, and $[V^{(2)}\backslash U_V,\partial(V^{(2)}\backslash U_V)]$ is
the fundamental class.
\end{lem}
\begin{proof}[Proof of lemma~\ref{lemma}]
Let us recall the construction of the class $(\widetilde{V}_2,f_2)$.
Since $f$ is an immersion, $\Delta(V)$ is a closed subset in
$(f^{(2)})^{-1}(\Delta(M))$. Since $f$ is a generic immersion,
$f^{(2)}$ is transversal to $\Delta(M)$ outside $\Delta(V)$. Therefore 
$(f^{(2)})^{-1}(\Delta(M))\backslash\Delta(V)$ is a compact
oriented submanifold without boundary
$\widetilde{f}_2:\widetilde{V}_2\hookrightarrow V^{(2)}\backslash\Delta(V)$.
Then the composition $\pi_1\circ\widetilde{f}_2$ is
$f_2:\widetilde{V}_2\to V^n$ (see details in \cite{Ronga}). By definition
\cite{Sw} of {\it Lefschetz duality}
$\gamma:\Omega^*(V^{(2)}\backslash U_V,\partial(V^{(2)}\backslash U_V))\to
\Omega_*(V^{(2)}\backslash U_V)$
\begin{equation*}
\begin{split}
(f^{(2)})^*((j^*)^{-1}t)\cap
\left[V^{(2)}\backslash U_V,\partial(V^{(2)}\backslash U_V)\right]=
(-1)^{2n\cdot m}&\gamma\left((f^{(2)})^*((j^*)^{-1}t)\right)\\&\qquad=
\gamma\left((f^{(2)})^*((j^*)^{-1}t)\right)
\end{split}
\end{equation*}
Since $f^{(2)}$ is transversal to $\Delta(M)$ outside $\Delta(V)$, we have
\begin{equation*}
\begin{split}
(\pi_1)_*\left(\gamma\left((f^{(2)})^*((j^*)^{-1}t)\right)\right)=
(\pi_1)_*&\left((\widetilde{f}_2)_*
\left[(f^{(2)})^{-1}(\Delta(M))\backslash\Delta(V)\right]\right)\\
&\qquad\qquad=(\pi_1)_*\left(\widetilde{V}_2,\widetilde{f}_2\right)=
(\widetilde{V}_2,f_2)
\end{split}
\end{equation*}
\end{proof}
\begin{proof}[Proof of lemma~\ref{double}]
To prove lemma~\ref{double}, it suffices to interpret the right hand side
of formula~(\ref{formula}) in terms of the differential $df$. Since
$\partial(V^{(2)}\backslash U_V)=\partial(U_V)$, we have
\begin{equation*}
\partial_*
\left[V^{(2)}\backslash U_V,\partial(V^{(2)}\backslash U_V)\right]=
\left[\partial(V^{(2)}\backslash U_V)\right],
\end{equation*}
where $\partial_*:\Omega_{2n}(V^{(2)}\backslash U_V,
\partial(V^{(2)}\backslash U_V))\to
\Omega_{2n-1}(\partial(V^{(2)}\backslash U_V))$ is the differential in
the exact bordism sequence of pair. Denote by $j_1$ the inclusion
$S\tau M\hookrightarrow\tau M\backslash\tau_0M$. Obviously, the map $j_1$
is a homotopy equivalence. Since
$j_1\circ s_M:M^m\to\tau M\backslash\tau_0M$ is a nowhere zero
cross-section of $\tau M$, we have
\begin{equation*}
\delta^*\left((j_1^*)^{-1}\gamma(M^m,s_M)\right)=t,
\end{equation*}
where $\delta^*:\Omega^{m-1}(\tau M\backslash\tau_0M)\to
\Omega^m(\tau M,\tau M\backslash\tau_0M)$ is the differential in the exact
cobordism sequence of pair, and  $\gamma$ is the Poincare duality
on the total manifold $S\tau M$. Denote by $j_2$ the isomorphism
$\partial(V^{(2)}\backslash U_V)\overset\sim\to S\tau V$.
From the explicit formula \cite{MS} for the ${\mathbb Z}_2$-equivariant
isomorphism, that identify a small neighborhood of zero section of
$\tau V$ with the neighborhood $U_V$
\begin{equation*}
\tau V\ni(x,\vec v_x)\mapsto(\exp_x(\vec v_x),\exp_x(-\vec v_x))\in V\times V,
\end{equation*}
it follows that the following diagram commutes (double arrows here denote
isomorphisms).
\begin{equation*}
\begin{matrix}
{\xymatrix{
{\Omega^{m-1}(S\tau V)}\ar@{=>}[r]^{j_2^*}&{\Omega^{m-1}(\partial(V^{(2)}\backslash U_V))}\ar[r]^{\delta^*}&{\Omega^m(V^{(2)}\backslash U_V,\partial(V^{(2)}\backslash U_V))}\\
&&{\Omega^m(M^{(2)},M^{(2)}\backslash\Delta(M))}\ar@{=>}[d]^{j^*}\ar[u]_{(f^{(2)})^*}\\
{\Omega^{m-1}(S\tau M)}\ar[uu]_{Sdf^*}&{\Omega^{m-1}(\tau M\backslash\tau_0M)}\ar[r]^{\delta^*}\ar@{=>}[l]_{j_1^*}&{\Omega^m(\tau M,\tau M\backslash\tau_0M)}\\
}}\\ \\
\text{Figure~1}
\end{matrix}
\end{equation*}
Therefore, $\delta^*(j_2^*Sdf^*\gamma(M^m,s_M))=(f^{(2)})^*((j^*)^{-1}t)$.
From the naturality of the $\cap$-product \cite{Sw} we have
\begin{equation*}
(j_3)_*\left(j_2^*Sdf^*\gamma(M^m,s_M)\cap
\left[\partial(V^{(2)}\backslash U_V)\right]\right)=
(f^{(2)})^*((j^*)^{-1}t)\cap
\left[V^{(2)}\backslash U_V,\partial(V^{(2)}\backslash U_V)\right],
\end{equation*}
where $j_3:\partial(V^{(2)}\backslash U_V)\to V^{(2)}\backslash U_V$
is the inclusion. Then, by lemma~\ref{lemma}
\begin{equation*}
\begin{split}
(\widetilde{V}_2,f_2)&=(\pi_1)_*\left((f^{(2)})^*((j^*)^{-1}t)\cap
\left[V^{(2)}\backslash U_V,\partial(V^{(2)}\backslash U_V)\right]\right)\\
&=(\pi_1\circ j_3)_*\left(j_2^*Sdf^*\gamma(M^m,s_M)\cap
\left[\partial(V^{(2)}\backslash U_V)\right]\right)\\
&=(-1)^{(2n-1)\cdot(m-1)}
(\pi_1\circ j_3)_*\gamma(j_2^*Sdf^*\gamma(M^m,s_M))\\
&=(-1)^{m-1}
(\pi_1\circ j_3)_*\circ (j_2)_*^{-1}\gamma Sdf^*\gamma(M^m,s_M)\\
&=(-1)^{m-1}
(\pi_1\circ j_3)_*\circ (j_2)_*^{-1}Sdf^!(M^m,s_M)
\end{split}
\end{equation*}
It remains only to note that $i_*=(\pi_1\circ j_3)_*\circ (j_2)_*^{-1}$.
\end{proof}

To prove the corollaries~\ref{ecor} and~\ref{maincor} we will need
the following results from the paper~\cite{Ronga}.
\begin{thm}[\cite{Ronga}]\label{herbert}
For any smooth generic immersion $f:V^n\looparrowright M^m$ of the
compact oriented manifold without boundary $V^n$ to the oriented
manifold without boundary $M^m$ such that $m-n$ is even, we have
\begin{equation}\label{herform}
v_k=f^*(m_{k-1})-e\cup v_{k-1},
\end{equation}
where $e$ is the Euler class of the normal bundle of the immersion
$f$ over $V^n$.
\end{thm}
\begin{cor}[\cite{Ronga}]\label{roncor}
Under the conditions of theorem~\ref{herbert}, we have
\begin{equation*}
(k-1)!\cdot v_k=\varphi_{k-1}\circ\varphi_{k-2}\circ\cdots\circ\varphi_1(1_V),
\end{equation*}
where $\varphi_k(a)=f^*f_!(a)-k\cdot e\cup a$; $e$ is the Euler class
of the normal bundle of the immersion $f$ over $V^n$.
\end{cor}
\begin{proof}[Proof of corollary~\ref{ecor}]
It suffices to substitute (\ref{center}) into formula~(\ref{herform})
with $k=2$.
\end{proof}
\begin{proof}[Proof of corollary~\ref{maincor}]
follows immediately from corollaries~\ref{roncor} and~\ref{ecor}.
\end{proof}

\section*{Acknowledgments}
I would like to thank P.M.~Akhmetiev and Yu.P.~Soloviev for useful remarks.

\end{document}